\documentclass[francais]{smfart}

\def\d{\succ}
\def\g{\prec}

\def\DD{\Delta}
\def\aa{\alpha}
\def\bb{\beta}
\def\cc{\gamma}
\def\dd{\delta}

\def\Prim{\operatorname{Prim}}
\def\Hom{\operatorname{Hom}}
\def\Id{\mathrm{Id}}

\def\ne{\nearrow}
\def\se{\searrow}
\def\sw{\swarrow}
\def\nw{\nwarrow}

\author{Jean-Louis Loday}
\address{Institut de Recherche Math\'ematique Avanc\'ee,
    CNRS et Universit\'e Louis Pasteur
    7 rue R. Descartes,
    67084 Strasbourg Cedex, France}
\email{loday@math.u-strasbg.fr}
\urladdr{www-irma.u-strasbg.fr/~loday/}
\title{Scindement d'associativit\'e et alg\`ebres de Hopf}
\alttitle{Splitting associativity and Hopf algebras} 
\subjclass{16A24, 16W30, 18D50}
\keywords{Alg\`ebre de Hopf, op\' erade, dendriforme, s\'erie g\'en\'eratrice}

\dedicatory{En hommage \`a Jean Leray}

\begin{abstract}
On montre que certaines alg\`ebres associatives dont le
produit se scinde en somme de plusieurs
op\'erations et qui sont libres, en un certain sens, pour ces op\'erations,
poss\`edent une structure d'alg\`ebre de Hopf. On
montre que l'op\'erade des  alg\`ebres dendriformes joue un r\^ole
particulier dans ce contexte, puis on donne de nombreux exemples.
\end{abstract}

\begin{altabstract}
  We show that some associative algebras whose product splits up into the
  sum of several operations and are free, in a certain sense, with respect
  to these operations, admit a Hopf algebra structure.  We show that the
  operad of dendriform algebras play a crucial role in this context, and we
  give numerous examples.
\end{altabstract}
\begin{document}
\frontmatter

\maketitle

\mainmatter

\section*{Introduction}

Dans leur c\'el\`ebre article sur les alg\`ebres de Hopf, John Milnor et John
Moore interpr\`ete le th\'eor\`eme 8 de l'article \cite{Leray} de Jean Leray de
la fa\c{c}on suivante (cf. th\'eor\`eme 7.5 de \cite{MM}): si une alg\`ebre
commutative unitaire $A$ poss\`ede une co-op\'eration unitaire, i.e. un
homomorphisme d'alg\`ebres associatives
\[
\Delta : A \to A \otimes A
\]
compatible avec l'unit\'e, alors $A$ est libre comme alg\`ebre associative
et commutative  (c'est-\`a-dire est une alg\`ebre
sym\'etrique). Ce r\'esultat peut s'\'etendre
\`a d'autres types d'alg\`ebres \`a condition de remplacer le produit
tensoriel par la somme (colimite) dans cette
cat\'egorie d'alg\`ebres (cf. Fresse \cite{Fr} et Oudom \cite{O}).

Le but de ce papier est, en un certain sens, de renverser la situation et
de montrer que, pour certains types d'alg\`ebres, on peut construire un
coproduit sur l'alg\`ebre libre. Dans le cas classique des alg\`ebres
associatives, l'alg\`ebre libre sur l'espace vectoriel $V$ est l'alg\`ebre
tensorielle $T(V)$ (alg\`ebre des polyn\^omes non commutatifs sur une base de
$V$). On sait que c'est aussi une alg\`ebre de Hopf pour le coproduit
construit \`a partir des shuffles. Comme cons\'equence importante de cette
propri\'et\'e les alg\`ebres enveloppantes des alg\`ebres de Lie sont des alg\`ebres
de Hopf. En pratique on peut utiliser le fait que $T(V)$ est libre pour
d\'emontrer la coassociativit\'e du coproduit shuffle sans calculs
combinatoires fastidieux. Nous montrons dans ce papier que cette technique
peut \^etre \'etendue \`a certains types d'alg\`ebres pr\'esentant un ``scindement
d'asso\-cia\-ti\-vi\-t\'e". Nous montrons que, lorsque certaines propri\'et\'es
de coh\'erence existent entre les relations d\'efinissant le type d'alg\`ebre et
l'unit\'e, alors l'alg\`ebre libre pour ce type (d\^ument augment\'ee) est une
alg\`ebre de Hopf. Dans le cas o\`u le type d'alg\`ebres est d\'efini par deux
op\'erations dont la somme est une op\'eration associative (scindement
d'associativit\'e), on constate que les relations doivent \^etre combinaisons
lin\'eaires de 3 relations particuli\`eres. Celles-ci sont exactement les
relations des ``alg\`ebres dendriformes".

Dans le premier paragraphe on explique ce qu'on entend par ``scindement
d'asso\-cia\-ti\-vi\-t\'e" et ``coh\'erence des relations avec l'unit\'e". On montre le
r\^ole primordial des alg\`ebres dendriformes pour ce probl\`eme. Dans le
deuxi\`eme paragraphe on d\'emontre l'existence d'une structure d'alg\`ebre de
Hopf sur les alg\`ebres libres pour les types d'alg\`ebres satisfaisant aux
conditions de coh\'erence.  Les deux premiers paragraphes sont restreints aux
types d'alg\`ebres ayant deux op\'erations g\'en\'eratrices sans sym\'etrie. On peut
\'etendre le r\'esultat \`a d'autres types d'alg\`ebres, ce qu'on fait dans le
troisi\`eme paragraphe, \'ecrit avec la terminologie des op\'erades qui est le
langage adapt\'e dans ce domaine.  L'existence du coproduit sur l'alg\`ebre
libre a pour application la g\'en\'eralisation de la notion de convolution.

Outre les alg\`ebres dendriformes, il se trouve que la plupart des nouveaux
types d'alg\`ebres avec scindement d'associativit\'e apparus derni\`erement
v\'erifient effectivement les propri\'et\'es de coh\'erence: alg\`ebres
2-associatives, trig\`ebres dendriformes, alg\`ebres pr\'e-dendriformes, alg\`ebres
dipt\`eres, quadrig\`ebres, alg\`ebres magmatiques. Apr\`es avoir donn\'e la
pr\'esentation de ces types d'alg\`ebres, on indique ce qui est connu sur leur
alg\`ebre libre dans le quatri\`eme paragraphe.

Dans le dernier paragraphe nous abordons le probl\`eme de la
d\'etermination de l'op\'erade des primitifs.

Je remercie Mar\' \i a Ronco et Teimuri Pirashvili pour les nombreuses
conversations et id\'ees \'echang\'ees sur le sujet, et Fran\c{c}ois Lamarche pour
une remarque pertinente.

\noindent \textbf{Notation.} Dans ce papier $K$ est un corps de caract\'eristique
quelconque. Par espace ou espace  vectoriel on entend un
espace vectoriel sur $K$. Le produit tensoriel sur $K$ des espaces $V$ et
$W$ est not\'e $V\otimes W$.

\section{Scindement d'associativit\'e et coh\'erence unitaire}
\label{sec:Scindement}

\subsection{D\'efinition}
\label{subsec:ScindDef}
 Soit $A$ une alg\`ebre associative (non
unitaire) dont on note $*$ le produit. On
dira qu'il y a \emph{scindement d'associativit\'e} lorsque cette op\'eration
$*$ est somme de deux op\'erations:
\begin{align*}
\label{formule0}
 x*y = x\g y + x \d y,  \tag{0} 
\end{align*}
que l'on qualifie respectivement de \emph{gauche} et \emph{droite}, et
lorsque l'associativit\'e de $*$ est
une cons\'equence des relations satisfaites par $\g$ et $\d$. L'exemple
suivant va jouer un r\^ole
primordial dans notre probl\'ematique.

\subsection{Exemple: alg\`ebres dendriformes}
\label{subsec:ScindEx}
Par d\'efinition une
\emph{alg\`ebre dendriforme} (cf. \cite{L2}), encore appel\'ee dig\`ebre
dendriforme, est un espace vectoriel $A$ muni de deux op\'erations
\emph{gauche} et \emph{droite} satisfaisant aux relations
\begin{align*}
(x \g y) \g z &= x \g (y * z),\tag{R1}\label{eq:R1}\\
(x \d y) \g z &= x \d (y \g z), \tag{R2}\label{eq:R2}\\
(x * y) \d z &= x \d (y \d z), \tag{R3}\label{eq:R3}
\end{align*}
Par addition des relations on constate que l'op\'eration $*$ est
associative, on a donc bien scindement d'associativit\'e.

\subsection{Compatibilit\'e entre relations et action de l'unit\'e}
\label{subsec:Compat}
Toute alg\`ebre associative $A$ peut \^etre rendue unitaire formellement en
posant $A_+= K. 1\oplus A$ (alg\`ebre augment\'ee) avec le produit associatif
induit par celui de $A$, $1$ \'etant l'unit\'e pour $*$. On se pose la question
de savoir si, lorsque le produit associatif est scind\'e, on peut \'etendre les
op\'erations $\g$ et $\d$ \`a tout $A_+$.  On doit avoir $a=1*a=1\g a + 1 \d a$
d'une part et $a=a*1=a\g 1 + a \d 1$ d'autre part. Faisons les choix
suivants pour l'action de 1 sur $a\in A$:
\[
 1\g a = 0,\quad 1 \d a = a,\quad  a\g 1 = a,\quad  a \d 1 = 0\
.\leqno (\dag)
\]
On ne peut pas \'etendre $\g$ et $\d$ \`a $K$ donc $1\g 1$ et $1\d 1$ ne sont
pas d\'efinis.  On voudrait que l'extension des op\'erations $\g$ et $\d$ \`a
l'alg\`ebre unitaire $A_+$ par les formules ci-dessus soit \emph{compatible},
i.e. que les relations satisfaites par $\g$ et $\d$ soient valables sur
$A_+$ pour autant que les termes soient d\'efinis. On dira alors que $A_+$
est une \emph{alg\`ebre augment\'ee}.

Dans un premier temps on suppose que les relations satisfaites par $\g$ et
$\d$ sont  quadratiques et r\'eguli\`eres
(voir paragraphe \ref{subsec:OperadeAlg}).    
Ceci signifie que les relations sont des combinaisons
lin\'eaires de mon\^omes du type $(x\circ_1 y)\circ_2 z$
et du type
$x\circ_1 (y\circ_2 z)$ o\`u  $\circ_1$ et $\circ_2$ sont soit $\g$ soit
$\d$ (il y a donc 8 mon\^omes
possibles). On remarquera que les alg\`ebres dendriformes sont de ce type.
Elles ont \'et\'e \'etudi\'ees dans \cite{L2}.

\subsection{Proposition}
\label{subsec:ScindProp}
\emph{L'extension des op\'erations $\g$ et $\d$
\`a  l'alg\`ebre unitaire $A_+$ est compatible si et seulement si les
relations
 satisfaites par $\g$ et $\d$ sont des combinaisons lin\'eaires des relations
\eqref{eq:R1}, \eqref{eq:R2} et \eqref{eq:R3} d\'ecrites ci-dessus.}
\begin{proof}[Preuve]
Soit
\begin{align*}
&\aa (x\g y)\g z + \bb (x\g y)\d z + \cc (x\d y)\g z + \dd (x\d y)\d
z= \\
&\qquad\qquad\qquad\aa' x\g (y\g z) + \bb' x\g (y\d z) + \cc' x\d (y\g z) + \dd' x\d
(y\d z)
\end{align*}
une relation o\`u $\aa, \bb$, etc, sont des scalaires. En rempla\c{c}ant $x$,
resp. $y$, resp. $z$ par $1$, on obtient
\begin{align*}
\cc = \cc', \quad \dd = \dd',\\
\aa = \bb', \quad \bb = \dd',\\
\aa = \aa', \quad \cc = \cc', 
\end{align*}
respectivement. En effet, par exemple pour $x=1$, on obtient
\[
\cc (b\g c) + \dd (b\d c) = \cc' (b\g c) + \dd' (b\d c)
\]
pour tous $b,c \in A_+$, d'o\`u  l'\'egalit\'e de la premi\`ere ligne.

On en d\'eduit que la relation de d\'epart est de la forme
\begin{align*}
\aa \big((x\g y)\g z - x\g (y\g z) - x\g (y\d z)\big) + \cc \big((x\d y)\g
z - x\d (y\g z) \big)\qquad\\
\qquad\qquad + \bb \big((x\g y)\d z +   (x\d y)\d z - x\d (y\d z)\big) = 0 
\end{align*}
c'est-\`a-dire une combinaison lin\'eaire des trois relations (Ri).
\end{proof}

On examine plusieurs types de dig\`ebres dans le second paragraphe.

\section{Structure d'alg\`ebre de Hopf sur les alg\`ebres libres}
\label{sec:Hopf}

On consid\`ere un type d'alg\`ebres $\mathcal{P}$ ayant deux op\'erations
g\'en\'eratrices $\g$ et $\d$
et dont les relations sont des combinaisons lin\'eaires de  \eqref{eq:R1}, \eqref{eq:R2} et
\eqref{eq:R3}. On suppose que l'on est en
pr\'esence d'un scindement d'associativit\'e, c'est-\`a-dire que
l'op\'eration $*$, d\'efinie par la formule $(0)$,  est associative.

Soient $A$ et $B$ deux  alg\`ebres de type $\mathcal{P}$, dont on note $A_+, B_+$
la \emph{$\mathcal{P}$-alg\`ebre augment\'ee}.

\subsection{Proposition (Coh\'erence)}
\label{subsec:Coherence}
\emph{Les formules ci-apr\`es font de $A_+\otimes B_+$
une $\mathcal{P}$-alg\`ebre augment\'ee (l'unit\'e \'etant $1\otimes 1$):
\begin{align*}
(a\otimes b) \circ (a'\otimes b') &:= (a*a')\otimes (b\circ b') 
\mbox{ \emph{si} }b\in B\mbox{ \emph{ou} }  b'\in B, \\
(a\otimes 1) \circ (a'\otimes 1) &:= (a\circ a')\otimes 1  , 
\end{align*}
o\`u $\circ= \g $ et $\d$ (ou une combinaison lin\'eaire quelconque
d'icelles), $a, a'\in A_+$ et $b,b'\in B_+$.}

On dit alors que le choix d'action de l'unit\'e est \emph{coh\'erent} avec
les relations (ainsi dans ce cas la compatibilit\'e implique
la coh\'erence). Il est pratique pour les calculs d'utiliser la  formule (abusive)  
\[
(a*a')\otimes (1\circ 1) = a\circ a'\otimes 1.
\]

\begin{proof}[Preuve]
On remarque que les formules impliquent imm\'ediatement
\[
(a\otimes b) * (a'\otimes b') = (a*a')\otimes (b* b') 
\]
dans tous les cas. Ainsi la structure d'alg\`ebre associative induite est
bien la structure
habituelle.

Soient $a,a',a'' \in A_+$ et $b,b',b''\in B_+$. Soit $i = 1 , 2$ ou $3$. On montre tout
d'abord que la relation 
(Ri) est v\'erifi\'ee pour $x=a\otimes b , y=a'\otimes b' ,z=a''\otimes b''$ dans les deux cas suivants

 -- l'un des \'el\'ements $b,b',b''$ vaut $1$ et les deux autres sont dans $B$,

 -- deux  des \'el\'ements $b,b',b''$ valent $1$ et le troisi\`eme est dans $B$.

Les parties gauche et droite de la relation \eqref{eq:R1} s'\'ecrivent respectivement
\[
 ((a\otimes b)\g (a'\otimes b'))\g (a''\otimes b'') = (a*a'*a'')\otimes ((b\g b')\g b'')
\]
et
\[
 (a\otimes b)\g ((a'\otimes b')*(a''\otimes b'')) = (a*a'*a'')\otimes (b\g (b'* b''))
\]
si $b\in B$ ou $b'\in B$. Si l'un des  $b,b',b''$
seulement vaut 1, alors  les composantes dans $B_+$  sont \'egales  par la
Proposition \ref{subsec:ScindProp}.  
Si  $b=1=b''$ et $b'\in B$ les deux termes sont nuls, et
 sont donc \'egaux. Si  $b'=1=b''$ et  $b\in B$ les deux termes valent $
(a*a'*a'')\otimes b$,  ils
 sont donc \'egaux.
Si  $b=1=b'$  et $b''\in B$ les parties gauche et droite s'\'ecrivent
respectivement
\[
 ((a\otimes 1)\g (a'\otimes 1))\g (a''\otimes 1) = ((a\g a')* a'')\otimes (1\g b'')= 0 
\]
et
\[
 (a\otimes 1)\g ((a'\otimes 1)*(a''\otimes 1)) = (a* a'*a'')\otimes  (1\g b'') =0 
\]
elles sont donc \'egales.

Les parties gauche et droite de la relation \eqref{eq:R2} s'\'ecrivent respectivement
\[
 ((a\otimes b)\d (a'\otimes b'))\g (a''\otimes b'') = (a*a'*a'')\otimes ((b\d b')\g b'')
\]
et
\[
 (a\otimes b)\d ((a'\otimes b')\g (a''\otimes b'')) = (a*a'*a'')\otimes (b\d (b'\g  b'')).
\]
Si l'un des  $b,b',b''$
seulement vaut 1, alors les composantes dans $B_+$  sont \'egales  par la
Proposition \ref{subsec:ScindProp}.  
Si  $b=1=b''$ et $b'\in B$ les deux termes
valent $ (a*a'* a'')\otimes b'$ et sont donc \'egaux.
Si  $b=1=b'$ et $b''\in B$ les parties gauche et droite s'\'ecrivent
respectivement
\[
 ((a\otimes 1)\d (a'\otimes 1))\g (a''\otimes b'') = ((a\d a')* a'')\otimes (1\g b'')=0
\]
et
\[
 (a\otimes 1)\d ((a'\otimes 1)\g (a''\otimes b'')) = (a\otimes 1)\d ((a'*a'')\otimes (1\g b''))=0
\]
elles sont donc \'egales.
Si  $b'=1=b''$ et $b\in B$ les parties gauche et droite s'\'ecrivent
respectivement
\[
 ((a\otimes b)\d (a'\otimes 1))\g (a''\otimes 1) = ((a*a')\otimes (b\d 1))\g (a''\otimes 1)=0 
\]
et
\begin{align*}
 (a\otimes b)\d ((a'\otimes 1)\g (a''\otimes 1)) = (a\otimes b )\d ((a'\g a''))\otimes 1)\qquad \\
\qquad\qquad\qquad = (a*(a'\g a''))\otimes (b\d 1)=0
\end{align*}
elles sont donc \'egales.

La v\'erification pour (R3) est en tous points analogue \`a celle de
\eqref{eq:R1}.

Supposons maintenant que $A$ et $B$ satisfont \`a la relation $\mathrm{(R)}:= \aa \mathrm{(R1)}
+ \bb \mathrm{(R2)} + \cc \mathrm{(R3)}$ pour des scalaires $\aa , \bb, \cc$. Si $b,b',b''$
sont dans $B$, alors (R) est v\'erifi\'ee pour $x=a\otimes b , y=a'\otimes b'
,z=a''\otimes b''$ car la relation (R) est valable dans $B$. Si au moins
l'un des $b,b',b''$ est dans $B$, alors (R) est v\'erifi\'ee car, d'apr\`es les
calculs pr\'ec\'edents, la relation (Ri) est v\'erifi\'ee pour tout $i$. Si
maintenant $b=b'=b''=1$, alors (R) est v\'erifi\'ee (toujours pour
$x=a\otimes b , y=a'\otimes b' ,z=a''\otimes b''$) car la relation (R)
est valable dans $A$.

En conclusion on a d\'emontr\'e que $A\otimes K \oplus K\otimes B \oplus A\otimes B$ est
une alg\`ebre de type
$\mathcal{P}$, et donc $A_+\otimes B_+$ est une $\mathcal{P}$-alg\`ebre
augment\'ee.
\end{proof}

\subsection{Alg\`ebre de Hopf connexe} On rappelle qu'une \emph{big\`ebre}
$\mathcal{H} = (\mathcal{H}, *, \Delta, u, c)$ est la donn\'ee
d'une structure d'alg\`ebre associative unitaire  $ (\mathcal{H}, *,  u)$,
d'une structure de cog\`ebre coassociative co-unitaire  $
(\mathcal{H}, \DD,  c)$ et on
suppose que
$\DD$ et $c$ sont  des homomorphismes d'alg\`ebres unitaires. On dit qu'une
big\`ebre est une \emph{alg\`ebre de Hopf} si elle poss\`ede une antipode.

Une big\`ebre $\mathcal{H}$ est dite \emph{connexe} si la filtration $F_r
\mathcal{H}$ est compl\`ete, i.e. si $\mathcal{H}= \bigcup_r F_r\mathcal{H}$,
pour $F_0\mathcal{H}:= K\cdot 1$ et
\[
 F_r\mathcal{H} := \{ x\in \mathcal{H} \vert \DD(x) - 1\otimes x - x \otimes 1 \in
F_{r-1}\mathcal{H}\otimes F_{r-1}\mathcal{H}\} .
\]
On montre ais\'ement qu'une big\`ebre connexe admet une antipode, donc il y
a \'equivalence entre big\`ebre connexe et
alg\`ebre de Hopf connexe.

\subsection{Th\'eor\`eme}
\label{subsec:ThmAlg}
 \emph{Soit  $\mathcal{P}$ un type d'alg\`ebres ayant
deux op\'erations g\'en\'eratrices $\g$ et $\d$
et dont les relations sont des combinaisons lin\'eaires de  \eqref{eq:R1}, \eqref{eq:R2} et
\eqref{eq:R3}, l'une d'elles \'etant \eqref{eq:R1}$+$\eqref{eq:R2}$+$\eqref{eq:R3} (scindement
d'associativit\'e). Alors l'alg\`ebre libre augment\'ee
${\mathcal{P}}(V)_+$ est munie naturellement d'une structure d'alg\`ebre de Hopf
connexe.}

\begin{proof}[Preuve]
Puisqu'il y a scindement d'associativit\'e, l'alg\`ebre  $\mathcal{P}(V)$ est une
alg\`ebre associative et $\mathcal{P}(V)_+$ est une alg\`ebre associative unitaire
augment\'ee. Il
nous faut donc construire la coop\'eration $\DD$ et d\'emontrer ses propri\'et\'es.

D'apr\`es la proposition \ref{subsec:Coherence} 
l'alg\`ebre associative $\mathcal{P}(V)_+\otimes \mathcal{P}(V)_+$ est munie
d'une structure de $\mathcal{P}$-alg\`ebre augment\'ee.

Consid\'erons maintenant l'application lin\'eaire
\[
\dd:V \to \mathcal{P}(V)_+\otimes \mathcal{P}(V)_+,\qquad  v \mapsto v\otimes 1 + 1 \otimes v .
\]
Il existe une et une seule extension de $\dd$ en un morphisme de
$\mathcal{P}$-alg\`ebres augment\'ees $ \Delta : \mathcal{P}(V)_+\to
\mathcal{P}(V)_+\otimes \mathcal{P}(V)_+$ car $\mathcal{P}(V)$ est la
$\mathcal{P}$-alg\`ebre libre sur $V$.

Puisque c'est un morphisme de $\mathcal{P}$-alg\`ebres augment\'ees, c'est, a
fortiori, un morphisme
d'alg\`ebres associatives augment\'ees.

Il nous reste \`a montrer que $\Delta $ est coassociatif. Les morphismes
$(\Delta \otimes \Id) \circ \Delta$ et $(\Id \otimes \Delta ) \circ \Delta$
\'etendent tous les deux l'application lin\'eaire $V\to
\mathcal{P}(V)_+{}^{\otimes 3}$ qui envoie $v$ sur $1\otimes 1\otimes v +1
\otimes v\otimes 1 +v \otimes 1\otimes 1$. Par unicit\'e de l'extension on
d\'eduit la coassociativit\'e de $\Delta$.

On a ainsi construit la structure de big\`ebre. Montrons maintenant que
cette big\`ebre est connexe. De par sa d\'efinition
l'alg\`ebre libre $\mathcal{P}(V)$ est de la forme $\mathcal{P}(V)= \oplus _{n\geq
1}\mathcal{P}(V)_n $ o\`u $\mathcal{P}(V)_n $ est l'espace engendr\'e
lin\'eairement par des produits quelconques de $n$ \'el\'ements de $V$. Posons $\overline \Delta (x) 
= \Delta(x) - x\otimes 1 - 1 \otimes x= \sum x_1\otimes x_2$. On a $\deg x = \deg x_1 + \deg x_2 $. Comme
 $\deg x_1 \geq 1$  et $\deg x_2 \geq 1$, on a $\deg x_1 < \deg x $ et $\deg x_2 < \deg x $. Ainsi, si $x\in 
 \mathcal{P}(V)_n$, on a ${\overline \Delta} ^n(x)=0$. Il s'en suit que $\mathcal{P}(V)_n \subset F_n \mathcal{P}(V)$. 
On a donc $\mathcal{P}(V)_+ =  \bigcup_r F_r{\mathcal{P}(V)_+}$ et $\mathcal{P}(V)_+$ est une
big\`ebre connexe. Donc c'est une alg\`ebre de Hopf connexe.
\end{proof}

\subsection{Remarques} En fait $\mathcal{P}(V)_+ $ a une structure plus fine que
simplement celle d'une alg\`ebre de Hopf, c'est une
\emph{$\mathcal{P}$-alg\`ebre de Hopf}. Ceci signifie que le coproduit $\Delta$ est
un morphisme de $\mathcal{P}$-alg\`ebres augment\'ees.

Puisque $\dd$ est sym\'etrique (i.e. $\dd= \tau \circ \dd$ o\`u $\tau (x\otimes
y) = y \otimes x$ ), on pourrait penser que $\DD$ est
co-commutative. Il n'en est rien car, bien que $\tau$ soit un automorphisme
d'alg\`ebre associative, ce n'est pas un
automorphisme d'alg\`ebre de type $\mathcal{P}$.

\subsection{Alg\`ebres dendriformes}
\label{subsec:AlgsEnd}
On en a donn\'e la pr\'esentation en
\ref{subsec:ScindEx}.   
Dans \cite{L2} on a montr\'e que l'alg\`ebre dendriforme libre sur un
g\'en\'erateur, not\'ee $Dend(K)$, a pour base lin\'eaire les arbres binaires
planaires. Donc la dimension de ses composantes homog\`enes est
\[
(1, 2, 5, 14, 42, 132,  \ldots, c_n, \ldots )
\]
o\`u $c_n = \frac{1}{n+1} \binom{2n}{n} $ est le nombre de Catalan (nombre d'arbres
binaires planaires \`a $n+1$ feuilles). 

En \ref{subsec:Compat} 
on a fait un choix pour l'action de 1. Si on identifie 1 \`a l'arbre sans
sommet interne, alors on constate que les formules de \cite{L2} s'\'etendent
sans obstruction.

 Dans \cite{LR1} on a
construit explicitement un coproduit
$\Delta '$ sur l'alg\`ebre associative augment\'ee $Dend(K)_+$  en
d\'ecrivant $\Delta '(t)$ pour tout arbre binaire
planaire $t$ par la  formule de r\'ecurrence
\[
 \Delta '(t\vee s) = \sum t_{(1)} * s_{(1)} \otimes  t_{(2)} \vee  s_{(2)} + t\vee
s \otimes 1.
\]
Ici  $\vee$ d\'esigne le greffage des arbres et on a adopt\'e la notation
de Sweedler $\Delta '(t) = \sum  t_{(1)} \otimes  t_{(2)}$.
Montrons que $\Delta '$ coincide avec la co-op\'eration  $\DD$ construite en
\ref{subsec:ThmAlg}. 
On rappelle que $t\vee s = t \d Y \g s $ o\`u $Y$ est le g\'en\'erateur de
$Dend(K)$. On a
\begin{align*}
\DD(t\vee s) &= \DD( t \d Y \g s)=  \DD(t) \d \DD(Y) \g\DD(s)\\
&= ( t_{(1)} \otimes t_{(2)}) \d ( 1\otimes Y + Y \otimes 1) \g  (s_{(1)} \otimes s_{(2)})\\
&= (t_{(1)}*1* s_{(1)}) \otimes (t_{(2)} \d Y \g  s_{(2)}) +  ( t_{(1)}*Y*
s_{(1)}) \otimes (t_{(2)} \d 1 \g  s_{(2)}) \\
&=  (t_{(1)}*s_{(1)}) \otimes (t_{(2)} \vee s_{(2)}) + t\vee s \otimes 1. 
\end{align*}
 On a ainsi montr\'e que $\Delta = \Delta '$. En application on obtient une d\'emonstration
plus simple de la co-associativit\'e de $\Delta '$.

Le coproduit des alg\`ebres dendriformes joue un r\^ole crucial dans le
travail de M. Ronco sur une
g\'en\'eralisation du th\'eor\`eme de Milnor-Moore (cf. \cite{R1} et \cite{R2}).

Une construction diff\'erente du produit  et du coproduit sur les arbres
planaires binaires a \'et\'e donn\'ee par Brouder et Frabetti dans \cite{BF}. Le
fait que ces deux
alg\`ebres sont isomorphes a \'et\'e d\'emontr\'e par  Holtkamp
\cite{H} et ind\'ependamment par Foissy \cite{Fo} (voir aussi \cite{AS}). On trouvera dans
\cite{A},  \cite{C} et \cite{E}  des liens entre les alg\`ebres
dendriformes et d'autres types de structure alg\'ebrique.

\subsection{Alg\`ebres dipt\`eres}  Par d\'efinition (cf. \cite{LR3}) les
\emph{alg\`ebres dipt\`eres}
 ont deux op\'erations g\'en\'eratrices
$*$ et $\d$ v\'erifiant les relations
\begin{align*}
(x*y)*z &= x*(y*z)\tag{as}\\
(x*y)\d z &= x\d (y\d z)\tag{dipt}
\end{align*}

Il est clair que l'on peut aussi les d\'efinir par les op\'erations
$\g$ et
$\d$ (via la formule \eqref{formule0}) 
et les relations \eqref{eq:R1}$+$\eqref{eq:R2} et \eqref{eq:R3}.

L'avantage de la premi\`ere
pr\'esentation est qu'elle a un sens aussi sur les ensembles. Dans \cite{LR3} on
a montr\'e que l'alg\`ebre dipt\`ere libre sur un
g\'en\'erateur not\'ee $Dipt(K)$  a pour base lin\'eaire deux copies des arbres
planaires, ses composantes homog\`enes ont donc pour
dimension
\[
(1, 2, 6, 22, 90, \ldots , 2C_{n-1}, \ldots )
\]
 o\`u  $C_n$ est le super-nombre de Catalan (nombre d'arbres planaires \`a $n+1$ feuilles). La
m\'ethode expos\'ee ci-dessus nous a permis de
construire facilement la structure d'alg\`ebre de Hopf. Cette structure
joue un r\^ole-cl\'e dans la d\'emonstration du
r\'esultat principal de \cite{LR3} qui est une g\'en\'eralisation du
th\'eor\`eme de Milnor-Moore au cas non-cocommutatif.

\subsection{Dig\`ebres sans nom}  Les op\'erations g\'en\'eratrices sont $\g$
et $\d$ et les relations sont
\eqref{eq:R1}$+$\eqref{eq:R3} et \eqref{eq:R2}. On peut aussi, \'evidemment remplacer $\mathrm{(R1)}+\mathrm{(R3)}$
par l'associativit\'e de $*$. On
trouve une structure diff\'erente de la pr\'ec\'edente bien que l'alg\`ebre
libre sur un g\'en\'erateur ait les m\^emes dimensions.

\subsection{Dig\`ebres associatives admissibles}  Les op\'erations
g\'en\'eratrices sont $\g$ et $\d$ et la
relation est \eqref{eq:R1}$+$\eqref{eq:R2}$+$\eqref{eq:R3}, c'est-\`a-dire l'associativit\'e de $*$.
L'alg\`ebre libre sur un
g\'en\'erateur a une base form\'ee des ``arbres binaires hybrides" (cf. \cite{Pa}). En effet, si on
admet que 2 est inversible, on peut changer de base d'op\'erations g\'en\'eratrices, et prendre
$*$ et $\cdot$ d\'efinie par    
$x\cdot y:= x\g y - x \d y$.  On trouve pour dimensions de l'objet libre sur un g\'en\'erateur
\[
(1, 2, 7, 31, 154, \ldots )
\]
En fait la s\'erie g\'en\'eratrice (cf. \ref{subsec:OperadeAlg})   
est l'inverse pour la composition de la s\'erie $f(x) = -1+x^2+\frac{1}{1+x}$.

On observe que l'op\'eration $\{x,y\}:= x\g y- y\d x$ munit l'espace
sous-jacent d'une structure d'\emph{alg\`ebre de Lie admissible} (cf.
\cite{Re}) car l'antisym\'etris\'e de $\{-,-\}$ est
\[
\{x,y\}-\{y,x\} = x\g y - y\d x -y\g x + x\d y = x*y - y*x = [x,y].
\]

\section{Op\'erades et alg\`ebres de Hopf}
\label{sec:Operades}

Dans la premi\`ere partie (paragraphes 1 et 2) on s'est restreint
volontairement \`a des types particuliers d'alg\`ebres: deux
op\'erations g\'en\'eratrices binaires et des relations
non-sym\'etriques (voir ci-dessous) avec scindement d'associativit\'e. Par
la m\^eme m\'ethode
on peut \'etendre le th\'eor\`eme \ref{subsec:ThmAlg} 
 \`a d'autres types d'alg\`ebres ayant
par exemple
plusieurs op\'erations g\'en\'eratrices et/ou des relations plus
g\'en\'erales. Afin de traiter le cas plus g\'en\'eral il est
pratique (voire n\'ecessaire) de se placer dans le cadre des op\'erades
alg\'ebriques. On rappelle tr\`es bri\`evement les
rudiments de cette th\'eorie en \ref{subsec:OperadeAlg}.

\subsection{Op\'erades alg\'ebriques}
\label{subsec:OperadeAlg}
Soit $\mathcal{P}$ un type d'alg\`ebres et
$\mathcal{P} (V)$ la $\mathcal{P}$-alg\`ebre libre sur l'espace vectoriel
$V$. On suppose que $\mathcal{P} (V)$ est de la forme $\mathcal{P} (V)=
\oplus_{n\geq 1}\mathcal{P} (n)\otimes _{S_n} V^{\otimes n}$ o\`u les
$\mathcal{P} (n)$ sont des $S_n$-modules \`a droite. Le groupe sym\'etrique
$S_n$ op\`ere \`a gauche sur $V^{\otimes n}$ par permutation des variables. On
consid\`ere $\mathcal{P}$ comme un endofoncteur de la cat\'egorie des espaces
vectoriels. La structure de $\mathcal{P}$-alg\`ebre libre de $\mathcal{P}(V)$
fournit une transformation de foncteurs $\cc : \mathcal{P}\circ \mathcal{P}
\to \mathcal{P}$ ainsi que $u: \Id \to \mathcal{P}$ v\'erifiant les axiomes
d'associativit\'e et d'unitalit\'e usuels.  Cette donn\'ee $(\mathcal{P}, \cc ,
u)$ est appel\'ee une \emph{op\'erade alg\'ebrique}. Une $\mathcal{P}$-alg\`ebre
est alors la donn\'ee d'un espace vectoriel $A$ et d'une application lin\'eaire
$\cc_A : \mathcal{P}(A) \to A $ telle que $\cc_A \circ \cc (A) = \cc_A
\circ \mathcal{P}(\cc_A)$ et $\cc_A \circ u(A) = \Id _A$.

L'espace $\mathcal{P}(n)$ est l'espace des op\'erations $n$-aires pour les
$\mathcal{P}$-alg\`ebres. On supposera ici qu'il n'y a qu'une (\`a
homoth\'etie pr\`es) op\'eration unaire, \`a savoir l'identit\'e: $\mathcal{P}(1)=
K.\Id$ .

On supposera aussi que toutes les op\'erations sont engendr\'ees (par
composition) par des op\'erations binaires, c'est-\`a-dire une famille de g\'en\'erateurs lin\'eaires
de $\mathcal{P}(2)$. On dira alors que l'op\'erade est \emph{binaire}. Si les relations
entre op\'erations sont cons\'equences de relations
ne faisant intervenir que des mon\^omes \`a 2 op\'erations, on dira que
l'op\'erade
est \emph{quadratique}.

Si les op\'erations binaires n'ont pas de sym\'etrie et que dans les
relations, les variables $x,y$
et $z$ apparaissent toujours dans le m\^eme ordre, on dira alors que
l'op\'erade est
\emph{r\'eguli\`ere}. L'espace $\mathcal{P}(n)$ est alors de la forme
$\mathcal{P}_n\otimes K[S_n]$ o\`u  $\mathcal{P}_n$ est un espace vectoriel et, en tant que
$S_n$-module,
$\mathcal{P}(n)$ est une somme de repr\'esentations r\'eguli\`eres. La famille des $\mathcal{P}_n$
munie de la composition est parfois appel\'ee une op\'erade non-sym\'etrique. Dans ce cas
l'alg\`ebre libre, et
donc l'op\'erade r\'eguli\`ere, sont enti\`erement d\'etermin\'ees par l'alg\`ebre libre
sur un g\'en\'erateur
$\mathcal{P}(K) = \oplus_{n\geq 1}\mathcal{P}_n$.

On dira qu'il y a \emph{scindement d'associativit\'e} si $\mathcal{P}(2)$
contient une op\'eration, not\'ee $(x,y) \mapsto x*y$, qui est associative.
Dans ce cas on peut toujours trouver une base de l'espace $\mathcal{P}(2)$
telle que la somme des vecteurs de bases soit l'op\'eration $*$. Les exemples
des paragraphes 1 et 2 sont des op\'erades binaires quadratiques
r\'eguli\`eres avec scindement d'associativit\'e.

La \emph{s\'erie g\'en\'eratrice} de l'op\'erade $\mathcal{P}$ est la fonction
\[
f^{\mathcal{P}}(x):= \sum (-1)^n\frac{\dim \, \mathcal{P}(n)}{n!} x^n.
\]
Si l'op\'erade est r\'eguli\`ere on a $f^{\mathcal{P}}(x):= \sum_{n\geq 1}
(-1)^n \dim \, \mathcal{P}_nx^n$. Dans les exemples de ce papier c'est la
s\'erie $(\dim \mathcal{P}_n)_{n\geq 1}$ que l'on donne.

Les op\'erades les plus int\'eressantes sont celles qui sont ``de Koszul". Une
condition n\'ecessaire pour la Koszulit\'e est que l'inverse, pour la
composition, de la s\'erie g\'en\'eratrice soit aussi une s\'erie \`a coefficients
entiers altern\'es.

\subsection{Actions de l'unit\'e, compatibilit\'e et coh\'erence} Soit $\mathcal{P}$ une op\'erade
binaire quadratique. Par \emph{action de l'unit\'e} on entend le choix de
deux applications lin\'eaires
\[
\alpha :\mathcal{P}(2) \to \mathcal{P}(1)= K , \qquad \beta :
\mathcal{P}(2) \to \mathcal{P}(1)= K ,
\]
qui permettent de donner un sens \`a $a\circ 1$ et $1\circ a$ respectivement pour toute
op\'eration $\circ \in \mathcal{P}(2)$ et tout $a\in
A$ o\`u $A$ est une $\mathcal{P}$-alg\`ebre:
\[
 a\circ 1 = \alpha (\circ) (a) , \qquad  1\circ a = \beta  (\circ) (a) .
\]
Lorsque l'op\'erade $\mathcal{P}$ est avec scindement d'associativit\'e on
suppose que l'on fait le choix $a*1=a=1*a$ pour $*$, c'est-\`a-dire
$\alpha(*) = \Id = \beta(*)$.

On dira que le choix d'action de l'unit\'e est \emph{compatible} avec les
relations de $\mathcal{P}$ si les relations sont encore valables sur
$A_+:=K.1\oplus A$ pour autant que les termes soient d\'efinis.

Consid\'erons l'espace $A\otimes K.1\oplus K.1\otimes B\oplus A\otimes
B$ o\`u $A$ et $B$ sont des $\mathcal{P}$-alg\`ebres. En utilisant le choix
d'action de l'unit\'e on \'etend les op\'erations binaires $\circ \in
\mathcal{P}(2)$ \`a cet espace en posant, comme en \ref{subsec:Coherence}: 
\begin{align*}
(a\otimes b) \circ (a'\otimes b') &:= (a*a')\otimes (b\circ b') \quad 
\mbox{ si }b\in B\mbox{ ou }b'\in B, \\
(a\otimes 1) \circ (a'\otimes 1) &:= (a\circ a')\otimes 1 \quad \mbox{sinon}.
\end{align*}
On dira que le choix d'action de l'unit\'e est \emph{coh\'erent} avec les
relations de $\mathcal{P}$ si  $A\otimes K.1\oplus K.1\otimes B\oplus A\otimes B$
muni de ces op\'erations est une $\mathcal{P}$-alg\`ebre.

Observons qu'une condition n\'ecessaire pour la coh\'erence est la
compatibilit\'e. Dans certains cas,
compatibilit\'e entra\^{\i}ne coh\'erence (cf. Proposition \ref{subsec:ScindProp})  
mais ce n'est pas toujours vrai.

Si $C_+$ est une autre $\mathcal{P}$-alg\`ebre augment\'ee, on v\'erifie que $(A_+\otimes B_+)\otimes C_+$ et
$A_+\otimes (B_+\otimes C_+)$ ont la m\^eme structure de $\mathcal{P}$-alg\`ebre augment\'ee.

\subsection{Th\'eor\`eme}
\label{subsec:ThmCoOp}
 \emph{Soit $\mathcal{P}$ une op\'erade binaire
quadratique non sym\'etrique. Toute action de
l'unit\'e coh\'erente avec   les relations de $\mathcal{P}$ permet de munir la
$\mathcal{P}$-alg\`ebre libre augment\'ee $\mathcal{P}(V)_+$ d'une
co-op\'eration co-associative (i.e. un coproduit)
\[
\Delta : \mathcal{P}(V)_+ \to \mathcal{P}(V)_+ \otimes \mathcal{P}(V)_+ 
\]
qui est un morphisme de $\mathcal{P}$-alg\`ebres augment\'ees.}

\emph{En particulier s'il y a scindement d'associativit\'e, $\mathcal{P}(V)_+$ est une
alg\`ebre de Hopf connexe.}

\begin{proof}[Preuve]
L'action de l'unit\'e permet, gr\^ace \`a l'hypoth\`ese de coh\'erence,
de munir  $\mathcal{P}(V)_+ \otimes \mathcal{P}(V)_+$ d'une structure de
$\mathcal{P}$-alg\`ebre augment\'ee. Le coproduit $\Delta$ est l'unique morphisme
de $\mathcal{P}$-alg\`ebres augment\'ees qui \'etend
l'application lin\'eaire
$v\mapsto 1\otimes v + v \otimes 1$. Le reste de la d\'emonstration est le m\^eme que
pour le th\'eor\`eme \ref{subsec:ThmAlg}. 
\end{proof}

\noindent\textbf{Remarque.} Sous l'hypoth\`ese ``r\'eguli\`ere" la premi\`ere formule
 de \ref{subsec:Coherence} 
 munit $A\otimes B$ d'une structure de $\mathcal{P}$-alg\`ebre. Sans cette
 hypoth\`ese, ce n'est plus automatique, mais il arrive encore parfois que ce
 soit vrai. Dans ce cas le th\'eor\`eme \ref{subsec:ThmCoOp} 
est encore valable (cf. exemple \ref{subsec:Zinbiel}). 

\subsection{Convolution op\'eradique} Soit $\mathcal{P}$ une op\'erade binaire quadratique r\'eguli\`ere 
munie d'une action coh\'erente de l'unit\'e. On va montrer que, comme dans le
cas associatif on peut munir l'espace des endomorphismes de la
$\mathcal{P}$-alg\`ebre libre de produits de convolution.

\subsection{Proposition} \emph{Pour tout espace vectoriel $V$ l'espace $\Hom_K(\mathcal{P} (V)_+, \mathcal{P} (V)_+)$ est muni d'une structure de $\mathcal{P}$-alg\`ebre unitaire.}

\begin{proof}[Preuve]
 Soit $\mu \in \mathcal{P} (2)$ une op\'eration g\'en\'eratrice. Pour tout couple d'applications lin\'eaires 
$f, g : \mathcal{P} (V)_+\to  \mathcal{P} (V)_+ $ on d\'efinit $\bar \mu (f,g) :  \mathcal{P} (V)_+\to  \mathcal{P} (V)_+ $ par la formule suivante:
\[
 \bar \mu (f,g):= \mu \circ (f\otimes g) \circ \Delta . 
\]
Il est imm\'ediat de v\'erifier que les relations satisfaites par les op\'erations g\'en\'eratrices $\mu$ de l'op\'erade
$\mathcal{P}$ sont aussi satisfaites par les op\'erations $\bar \mu$. Donc 
 $\Hom_K(\mathcal{P} (V)_+, \mathcal{P} (V)_+)$ est une
 $\mathcal{P}$-alg\`ebre unitaire.
\end{proof}

\subsection{Remarque} Dans le cas des alg\`ebres associatives, i.e.~$\mathcal{P} = \mathbf{As}$ et $\mu = *$, l'op\'eration 
$\bar \mu$ est la convolution classique.

\section{Exemples}
\label{sec:Exemples}
\subsection{Alg\`ebres de Zinbiel}
\label{subsec:Zinbiel}
 Supposons donn\'ee une seule op\'eration
$\d$. On choisit pour
action de 1 les formules suivantes:
\[
1\d x = x , \qquad x\d 1 = 0 .
\]
On montre, comme dans la proposition \ref{subsec:ScindProp}, 
que l'unique relation possible
pour avoir coh\'erence
est
\[
(x\d y + y \d x) \d z = x \d ( y\d z ).
\]
Les alg\`ebres d\'efinies par $\d$ et cette relation sont les \emph{alg\`ebres de Zinbiel} (cf. \cite{L1}), qui sont
duales, au sens op\'eradique, des alg\`ebres de Leibniz. On constate que
l'op\'eration $*$ d\'efinie par
$x*y= x\d y + y \d x$ est associative. On peut montrer (cf. loc.cit.) que
l'alg\`ebre de Zinbiel libre augment\'ee
est l'alg\`ebre des shuffles $T^{sh}(V)$ (i.e. $T(V)$ en tant qu'espace
vectoriel, mais avec le
produit shuffle). La structure d'alg\`ebre de Hopf donn\'ee par notre
m\'ethode n'est rien d'autre que
la structure connue: $\DD$ est la d\'econcat\'enation. La preuve consiste \`a identifier 
$x_1x_2x_3  \ldots x_n$ \`a $(\ldots ((x_1\d x_2)\d x_3)  \ldots x_n)$ et raisonner par r\'ecurrence.

\subsection{Alg\`ebres pr\'e-dendriformes} Par d\'efinition une \emph{alg\`ebre
pr\'e-dendriforme} a 3 op\'erations $\g, \d, *$ qui
satisfont \`a 4 relations:
\begin{align*}
(x * y) *z &= x * (y * z),\tag{R0} \label{eq:R0} \\
(x \g y) \g z &= x \g (y * z),\tag{\ref{eq:R1}} \\
(x \d y) \g z &= x \d (y \g z),\tag{\ref{eq:R2}}\\
(x * y) \d z &= x \d (y \d z).\tag{\ref{eq:R3}}
\end{align*}

L'op\'erade associ\'ee est donc binaire quadratique r\'eguli\`ere avec
scindement d'associativit\'e.
On obtient les alg\`ebres dendriformes comme quotient en introduisant la
relation de sym\'etrie
\[
x*y = x\g y + x\d y,
\]
car sous cette condition la relation \eqref{eq:R0} devient \'egale \`a
\eqref{eq:R1}$+$\eqref{eq:R2}$+$\eqref{eq:R3}.

L'op\'erade des alg\`ebres pr\'edendriformes a pour dimension des parties
homog\`enes
\[
(1, 3, 14, 80,  510, \ldots ),
\]
dont la s\'erie g\'en\'eratrice altern\'ee est l'inverse pour la composition de
\[
-x+3x^2-4x^3+5x^4-\ldots = -1-x+\frac{1}{(1+x)^2}.
\]
Si on prend les m\^emes
conventions que dans le paragraphe 1, \`a savoir
\[
1\g a = 0,1\d a = a,a \g 1 = a ,a \d 1 = 0 ,1*a=a=a*1,
\]
on constate imm\'ediatement que les quatre relations \eqref{eq:R0} \`a \eqref{eq:R3} sont
compatibles. En fait ce choix est m\^eme
coh\'erent avec les relations et donc l'alg\`ebre pr\'e-dendriforme libre
augment\'ee,
 ${preDend}(V)_+$ est munie d'une structure d'alg\`ebre de Hopf par le
th\'eor\`eme \ref{subsec:ThmCoOp}. 

La version cog\`ebre des relations ci-dessus joue un r\^ole primordial dans
le travail de P. Leroux \cite{Le1, Le2}.

\subsection{Remarque \textrm{(due \`a F. Lamarche)}} Les formules \eqref{eq:R0} \`a
\eqref{eq:R3} se rencontrent
dans le contexte des foncteurs adjoints de la mani\`ere suivante.
Soit $(\mathcal{C}, *, I)$ une cat\'egorie mono\"{\i}dale de produit associatif $*$
et d'unit\'e $I$. Supposons que, pour tout objet $C$ de
$\mathcal{C}$,  le foncteur $C* - : \mathcal{C}\to \mathcal{C}$ ait un adjoint \`a
droite, que l'on note $-\g C$, et que le foncteur
$-* C : \mathcal{C}\to \mathcal{C}$ ait aussi un adjoint \`a droite  que l'on note $C\d
-$. On a donc:
\[
\Hom_\mathcal{C}(C*A,X)= \Hom_\mathcal{C}(A,X\g C), \quad  \Hom_{\mathcal{C}}(A*C,X)= \Hom_{\mathcal{C}}(A,C\d X).
\]
Alors, en \'evaluant de plusieurs
mani\`eres l'ensemble $\Hom_\mathcal{C}(A*B*C,X)$, on trouve pr\'ecis\'ement
les relations \eqref{eq:R1}, \eqref{eq:R2} et
\eqref{eq:R3}. Par exemple on a d'une part
\[
\Hom_\mathcal{C}(A*B*C,X)= \Hom_{\mathcal{C}}(A, (B*C)\d X)
\]
 et d'autre part
\[
\Hom_\mathcal{C}(A*B*C,X)= \Hom_\mathcal{C}(A*B, C\d X)= \Hom_\mathcal{C}(A, B\d
(C\d X)),
\]
ce qui donne la
relation \eqref{eq:R3}.

\subsection{Trig\`ebres dendriformes  \cite{LR2}} On se donne 3 op\'erations $\g,
\d $ et $\cdot$ et 7 relations:
\begin{align*}
&\begin{cases}
(x \g y) \g z = x \g (y * z), \\
(x \d y) \g z = x \d (y \g z), \\
(x * y) \d z = x \d (y \d z), \\
\end{cases}\\
&\begin{cases}
(x \d y) \cdot z = x \d (y \cdot z), \\
(x \g y) \cdot z = x \cdot (y \d z), \\
(x \cdot y) \g z = x \cdot (y \g z), 
\end{cases}\\
&\hspace{0.3cm}(x \cdot y) \cdot z = x \cdot (y \cdot z),
\end{align*}
o\`u $x*y := x\g y + x\d y  + x\cdot y$.

L'op\'erade des trig\`ebres dendriformes, not\'ee $Tridend$, est binaire, quadratique, r\'eguli\`ere, avec scindement
d'associativit\'e. 

On \'etend ces op\'erations \`a l'unit\'e par les choix suivants
\[
1\g a = 0, 1\d a = a, a \g 1 = a , a \d 1 = 0 , 1\cdot a=0=a\cdot
1.
\]
Il s'ensuit que l'on a bien $1*a=a=a*1$.

On peut montrer que ces choix sont coh\'erents avec les
relations. Ainsi la trig\`ebre
dendriforme libre unitaire $Tridend(V)_+$ peut \^etre munie d'une structure
d'alg\`ebre de Hopf.
Rappelons que $Tridend(V)$ se d\'ecrit explicitement \`a l'aide des arbres
planaires (cf. \cite{LR2}). La dimension  de ses parties homog\`enes est donc
donn\'ee par le super-nombre de Catalan $C_n$:
\[
(1, 3, 11, 45, 197, \ldots, C_{n-1}, \ldots  ).
\]
On peut expliciter  $\Delta$  sur les arbres planaires comme en
\ref{subsec:AlgsEnd} 
en utilisant les formules
\begin{align*}
x^0\vee \ldots \vee x^k &= ( x^0\vee \ldots \vee x^{k-1})\cdot (Y\g x^k), \hbox {
si }k>1 \hbox{  ou } k=1 \hbox { et } x^0\neq \vert,\\
\vert \vee x &= Y \g x .
\end{align*}
On constate que les \'el\'ements qui s'\'ecrivent
$\omega\cdot \theta$ avec $\omega$ et $\theta$
\'el\'ements primitifs de $\mathcal{P}(V)$ sont aussi des \'el\'ements primitifs.

\subsection{Alg\`ebres 2-associatives \cite{LR3} \cite{Pi}}  On se donne 2 op\'erations
associatives $*$ et $\cdot$ et pas d'autres relations. L'op\'erade des alg\`ebres
2-associatives, not\'ee $2as$, est binaire, quadratique, r\'eguli\`ere.
L'alg\`ebre libre sur un g\'en\'erateur est, comme pour les alg\`ebres dipt\`eres, de
dimension $2C_{n-1}$ en dimension $n\geq 2$. Ici on va modifier quelque peu
notre construction de $\Delta$.  On fait le choix d'actions de 1 suivant:
$1*a=a=a*1$ et $1\cdot a=a=a\cdot 1$ et on met sur le produit tensoriel de
deux alg\`ebres les produits diagonaux classiques. Il y a alors une et une
seule co-op\'eration unitaire
\[
\Delta : 2as(V)_+ \to 2as(V)_+\otimes 2as(V)_+
\]
qui v\'erifie
\[
  \begin{cases}
\DD(x*y) = \DD(x) *\DD(y), \\
\DD(x\cdot y) = (x\otimes 1)\cdot \DD( y) + \DD(x)\cdot  (1\otimes y) - x\otimes y,
  \end{cases}
\]
et $\DD(v) = 1\otimes v + v \otimes 1$ pour $v\in V$.
On v\'erifie que cette co-op\'eration est co-associative, co-unitaire, et est un morphisme pour $*$,
donc  $(2as(V), *, \DD)$ est une alg\`ebre de Hopf. Sa partie primitive est \'etudi\'ee dans \cite{LR3}.

La relation satisfaite entre la co-op\'eration $\DD$ et l'op\'eration $\cdot$ est
appel\'ee \emph{relation de Hopf infinit\'esimale unitaire}. La situation
typique est le cas du module tensoriel $T(V)$ \'equip\'e de $\cdot =$
concat\'enation et de $\DD=$ d\'econcat\'enation.

Il y a de nombreux exemples int\'eressants de quotient de cette op\'erade,
c'est-\`a-dire des op\'erades obtenues en rajoutant des relations. Par exemple
si on rajoute la relation
\[
(x * y) \cdot z = x *(y \cdot z),
\]
on obtient une op\'erade tr\`es similaire \`a l'op\'erade des alg\`ebres
dendriformes. En fait ces deux op\'erades sont reli\'ees par une homotopie.
Donc l'alg\`ebre libre sur un g\'en\'erateur est aussi index\'ee par les arbres
binaires planaires (cf. \cite{Pi}).  Cette op\'erade est binaire,
quadratique, r\'eguli\`ere. On obtient une structure d'alg\`ebre de Hopf sur
l'alg\`ebre libre.

\subsection{Quadrig\`ebres \cite{AL}} On se  donne 4 op\'erations  $\se$, $\ne$,
$\nw$ et$\sw$ et on note
\begin{align*}
x\d y:= x\ne y+x\se y,\\
x\g y:= x\nw y+x\sw y,\\
x\vee y:= x\se y+x\sw y,\\
x\wedge y:= x\ne y+x\nw y .  
\end{align*}
ainsi que
\[
x* y:=x\se y+x\ne y+x\nw y+x\sw y,
=   x\d y+x\g y=x\vee y+x\wedge y.
\]
Puis on se donne 9 relations:
\begin{align*}
(x\nw y)\nw z =x\nw(y* z),&  & (x\ne y)\nw z = x\ne(y\g z),  && (x\wedge y)\ne z=x\ne(y\d z), \\
 (x\sw y)\nw z =x\sw(y\wedge  z), & & (x\se y)\nw z = x\se(y\nw z),  && (x\vee y)\ne z=x\se(y\ne z),\\
(x\g y)\sw z=x\sw(y\vee  z), & & (x\d y)\sw z = x\se(y\sw z),  && (x* y)\se z =x\se(y\se z). 
\end{align*}
 L'op\'erade des quadrig\`ebres, not\'ee $\mathcal{Q}$, est binaire,
quadratique et r\'eguli\`ere.  On conjecture que sa s\'erie g\'en\'eratrice 
est l'inverse de la s\'erie $\sum_{n\geq 1} (-1)^n
n^2x^n= \frac{x(-1+x)}{(1+x)^3}$  pour la composition, ce qui donnerait pour
dimensions des composantes homog\`enes de la quadrig\`ebre libre sur un g\'en\'erateur:
\[
(1, 4, 23, 156, 1162, \ldots )
\]
On \'etend les quatre produits \`a l'unit\'e avec les choix suivants:
\[
 a\nw 1 = a, \quad 1\se a = a
\]
et tous les autres produits avec 1 sont nuls. Ainsi on a
\[
a\g 1 =a,1\d a = a,a\wedge 1 = a,1\vee a =a ,a*1=a=1*a,
\]
et les autres produits  nuls.
On v\'erifie imm\'ediatement que les 9 relations sont compatibles avec ces
choix. On peut montrer qu'elles sont m\^eme coh\'erentes et on   en d\'eduit
que l'alg\`ebre libre augment\'ee $\mathcal{Q}(V)_+$ est une alg\`ebre de Hopf.

Lorsque les 4 op\'erations g\'en\'eratrices satisfont aux propri\'et\'es de
sym\'etrie
\[
x\nw y = y\se x, \quad x\sw y = y \ne x , 
\]
on dit que les quadrig\`ebres sont \emph{commutatives}. Nos choix sont
compatibles avec ces
relations de sym\'etrie car $a\nw 1 = a = 1\se a$ et tous les autres termes
contenant $1$ sont nuls. On a donc une structure
d'alg\`ebre de Hopf sur la quadrig\`ebre commutative libre augment\'ee.

\subsection{Alg\`ebres magmatiques \cite{GH}}  Donnons-nous une op\'eration
binaire, not\'ee $(x,y) \mapsto x\cdot y$, sans aucune
relation. L'alg\`ebre libre sur un g\'en\'erateur admet \'evidemment les arbres binaires
planaires pour base lin\'eaire. Elle a donc pour dimensions
(cf. \ref{subsec:AlgsEnd}): 
\[
(1,1,2,5,14,\ldots , c_{n-1}, \ldots ).
\]
Ce cas est un peu diff\'erent
des pr\'ec\'edents puisqu'on n'a plus de scindement d'asso\-ciati\-vit\'e. Prenons le choix d'action de
l'unit\'e usuel:
$1\cdot a = a = a\cdot 1$. On v\'erifie imm\'ediatement que les conditions de coh\'erence
sont v\'erifi\'ees et on peut appliquer le th\'eor\`eme \ref{subsec:ThmCoOp}. 
L'int\'er\^et de ce cas r\'eside dans la nature des \'el\'ements primitifs d\'efinis
par le coproduit (voir paragraphe suivant).

Une variante de l'op\'erade magmatique est l'op\'erade \emph{magmatique
commutative} o\`u l'on suppose de plus que $a\cdot b =
b \cdot a$ pour tous \'el\'ements $a$ et $b$. On a encore coh\'erence dans
ce cas.

Un autre quotient int\'eressant de l'op\'erade magmatique est l'op\'erade des alg\`ebres pre-Lie
 (cf.~\cite{CL}), qui sont caract\'eris\'ees par la relation:
\[
(x\cdot y)\cdot z - x\cdot (y \cdot z ) = (x\cdot z)\cdot y - x\cdot (z \cdot y ) .
\]
La coh\'erence est valide dans ce cas bien que l'op\'erade ne soit pas r\'eguli\`ere.

\subsection{Op\'erades ensemblistes et arithm\'etique} Lorsque l'op\'erade alg\'ebrique provient
d'une op\'erade ensembliste, c'est-\`a-dire bien d\'efinie sur la cat\'egorie des
ensembles, on peut construire une arithm\'etique sur l'objet ensembliste
libre. Observons que plusieurs des op\'erades pr\'esent\'ees dans les sections
pr\'ec\'edentes sont ensemblistes (alg\`ebres dipt\`eres, alg\`ebres
pr\'e-dendriformes, alg\`ebres 2-associatives, alg\`ebres magmatiques). On se
sert de l'op\'eration associative pour construire l'addition, et de la
composition dans l'alg\`ebre libre pour construire la multiplication. Pour
l'op\'erade $As$ c'est tout simplement l'arithm\'etique sur $\mathbb{N}$. Pour
l'op\'erade magmatique (cf. 4.7) c'est l'arithm\'etique sur les arbres
planaires binaires \'evoqu\'ee dans \cite{B}. M\^eme quand l'op\'erade n'est pas
ensembliste on peut parfois construire une arithm\'etique gr\^ace \`a un bon
choix de base lin\'eaire de l'alg\`ebre libre. C'est ce qui est fait dans
\cite{L3} pour les dig\`ebres dendriformes et les trig\`ebres dendriformes.

\section{L'op\'erade des primitifs}
\label{sec:primitifs}

Soit $\mathcal{P}$ une op\'erade binaire
quadratique pour laquelle, apr\`es avoir fait un choix
d'action de l'unit\'e, on a r\'eussi \`a construire une co-op\'eration
co-associative sur l'alg\`ebre libre augment\'ee:
\[
\Delta : \mathcal{P}(V)_+ \to \mathcal{P}(V)_+\otimes \mathcal{P}(V)_+
\]
On peut alors d\'efinir l'espace des \'el\'ements primitifs par
\[
\Prim \mathcal{P}(V) = \{ x \in \mathcal{P}(V) \vert \DD(x) = 1\otimes x + x \otimes 1\}.
\]
Lorsque $\Delta$ est un morphisme de $ \mathcal{P}$-alg\`ebres augment\'ees, on peut montrer que la composition
d'\'el\'ements primitifs est encore un \'el\'ement primitif. Donc $\Prim \mathcal{P}(V)$ est l'alg\`ebre libre d'un certaine op\'erade $\Prim \mathcal{P}$.

Le jeu consiste maintenant \`a trouver quelle est cette op\'erade dans les cas
qui nous int\'eressent. L'outil principal est l'\emph{idempotent Eul\'erien} dans le cas 
cocommutatif et  l'\emph{idempotent de Ronco} dans le cas non-cocommutatif 
( cf. \cite{R1} et \cite{R2}). Voici quelques r\'eponses:
\begin{align*}
\Prim As &= Lie \\
 \Prim  Com  &= Vect  \\
 \Prim  Dend &=   \hbox{op\'erade des  alg\`ebres braces, cf. \cite{R1},} \\
 \Prim Dipt  &=  \hbox{op\'erade des $B_{\infty}$-alg\`ebres, cf. \cite{LR3},}\\
 \Prim 2as  &=  \hbox{op\'erade des $B_{\infty}$-alg\`ebres, cf. \cite{LR3}.}
\end{align*}
Dans le cas des alg\`ebres magmatiques (cf.~4.7), les premiers calculs ont
\'et\'e faits par Gerritzen et Holtkamp \cite{GH}.
L'op\'erade $\Prim Mag$ dont l'alg\`ebre libre est la partie primitive de
l'alg\`ebre magmatique libre poss\`ede au moins une
op\'eration binaire antisym\'etrique, notons la $[-,-]$, et une
op\'eration ternaire, notons-la $as(-,-,-)$ car si $x,y$ et $z$ sont
primitifs, alors il en est de m\^eme de $[x,y] := x\cdot y - y \cdot x$ et
de $as(x,y,z) := (x\cdot y) \cdot z - x \cdot (y\cdot z)$.
Ces deux op\'erations g\'en\'eratrices sont li\'ees par la relation
\begin{align*}
as(x,y,z)+ as(y,z,x)+ as(z,x,y) - as(x,z,y)- as(y,x,z)- as(z,y,x) = \qquad\\
\qquad\qquad  [[x,y],z]+[[y,z],x]+[[z,x],y]
\end{align*}
qu'on pourrait appeler la \emph{relation de Jacobi non-associative}.
Mais il est montr\'e dans \cite{GH} que ce n'est pas suffisant, et qu'il y a
d'autres g\'en\'erateurs, par exemple l'\'el\'ement de l'alg\`ebre
magmatique libre
\[
as(x, y, z\cdot t)  - z\cdot as(x,y,t)- as(x,y,z)\cdot t
\]
est primitif, mais n'est pas engendr\'e par les op\'erations crochet et
associateur.

\backmatter

\bibliographystyle{smfalpha}
\bibliography{jllodayScind}

\providecommand{\bysame}{\leavevmode ---\ }
\providecommand{\og}{``}
\providecommand{\fg}{''}
\providecommand{\smfandname}{\&}
\providecommand{\smfedsname}{\'eds.}
\providecommand{\smfedname}{\'ed.}
\providecommand{\smfmastersthesisname}{M\'emoire}
\providecommand{\smfphdthesisname}{Th\`ese}
\begin{thebibliography}{Rem02}

\bibitem[Agu00]{A}
{\scshape M.~Aguiar} -- {\og Pre-{P}oisson algebras\fg}, \emph{Letter Math.
  Phys.} \textbf{54} (2000), p.~263--277.

\bibitem[AL04]{AL}
{\scshape M.~Aguiar {\normalfont \smfandname} J.-L. Loday} -- {\og
  Quadri-algebras\fg}, to appear, 2004.

\bibitem[AS]{AS}
{\scshape M.~Aguiar {\normalfont \smfandname} F.~Sottile} -- {\og The
  {L}oday-{R}onco {H}opf algebra of planar binary trees\fg}, in progress.

\bibitem[BF03]{BF}
{\scshape C.~Brouder {\normalfont \smfandname} A.~Frabetti} -- {\og {QED}
  {H}opf algebras on planar binary trees\fg}, \emph{J. of Algebra} (2003),
  p.~298--322.

\bibitem[Blo95]{B}
{\scshape V.~Blondel} -- {\og Une famille d'op\'erations sur les arbres
  binaires\fg}, \emph{C. R. Acad. Sci. Paris} \textbf{321} (1995), p.~491--494.

\bibitem[Cha02]{C}
{\scshape F.~Chapoton} -- {\og Construction de certaines op\'erades et
  big\`ebres associ\'ees aux polytopes de {S}tasheff et hypercubes\fg},
  \emph{Trans. Amer. Math. Soc.} \textbf{354} (2002), p.~63--74.

\bibitem[CL01]{CL}
{\scshape F.~Chapoton {\normalfont \smfandname} M.~Livernet} -- {\og Pre-{L}ie
  algebras and the rooted trees operad\fg}, \emph{Int. Math. Res. Not.}
  \textbf{8} (2001), p.~395--408.

\bibitem[EF02]{E}
{\scshape K.~Ebrahimi-Fard} -- {\og Loday-type algebras and the {R}ota-{B}axter
  relations\fg}, \emph{Letter Math. Phys.} \textbf{61} (2002), p.~139--147.

\bibitem[Foi02]{Fo}
{\scshape L.~Foissy} -- {\og Les alg\`ebres de {H}opf des arbres enracin\'es
  d\'ecor\'es. {II}\fg}, \emph{Bull. Sci. Math.} \textbf{126} (2002),
  p.~249--288.

\bibitem[Fre98]{Fr}
{\scshape B.~Fresse} -- {\og Cogroups in algebras over an operad are free
  algebras\fg}, \emph{Comment. Math. Helv.} \textbf{73} (1998), no.~4,
  p.~637--676.

\bibitem[GH03]{GH}
{\scshape L.~Gerritzen {\normalfont \smfandname} R.~Holtkamp} -- {\og Hopf
  co-addition for free magma algebras and the non-associative {H}ausdorff
  series\fg}, \emph{J. of Algebra} \textbf{265} (2003), p.~264--284.

\bibitem[Hol03]{H}
{\scshape R.~Holtkamp} -- {\og Comparison of {H}opf algebras on trees\fg},
  \emph{Arch. Math. (Basel)} \textbf{80} (2003), p.~368--383.

\bibitem[Ler45]{Leray}
{\scshape J.~Leray} -- {\og Sur la forme des espaces topologiques et sur les
  points fixes des repr\'esentations\fg}, \emph{J. Math. Pures App.}
  \textbf{24} (1945), p.~95--167.

\bibitem[Ler02]{Le1}
{\scshape P.~Leroux} -- {\og Coassociativity breaking and oriented graphs\fg},
  preprint, ArXiv: {\tt math.QA/0204342}, 2002.

\bibitem[Ler03]{Le2}
\bysame , {\og From entangled codipterous coalgebras to coassociative
  manifolds\fg}, preprint, ArXiv: {\tt math.QA/0301080}, 2003.

\bibitem[Lod95]{L1}
{\scshape J.-L. Loday} -- {\og Cup-product for {L}eibniz cohomology and dual
  {L}eibniz algebras\fg}, \emph{Math. Scand.} \textbf{77} (1995), no.~2,
  p.~189--196.

\bibitem[Lod01]{L2}
\bysame , {\og Dialgebras\fg}, in \emph{Dialgebras and related operads},
  Lecture Notes in Math., no. 1763, Springer, Berlin, 2001, p.~7--66.

\bibitem[Lod02]{L3}
\bysame , {\og Arithmetree\fg}, \emph{J. of Algebra} \textbf{258} (2002),
  no.~1, p.~275--309.

\bibitem[LR98]{LR1}
{\scshape J.-L. Loday {\normalfont \smfandname} M.~Ronco} -- {\og Hopf algebra
  of the planar binary trees\fg}, \emph{Adv. in Maths} \textbf{139} (1998),
  p.~293--309.

\bibitem[LR02]{LR2}
\bysame , {\og Trialgebras and families of polytopes\fg}, preprint, ArXiv: {\tt
  math.AT/0205043}, 2002.

\bibitem[LR03]{LR3}
\bysame , {\og Alg\`ebres de {H}opf colibres\fg}, \emph{C. R. Acad. Sci. Paris}
  \textbf{337} (2003), p.~153--158.

\bibitem[MM65]{MM}
{\scshape J.~W. Milnor {\normalfont \smfandname} J.~C. Moore} -- {\og On the
  structure of {H}opf algebras\fg}, \emph{Ann. of Math.} \textbf{81} (1965),
  p.~211--264.

\bibitem[Oud99]{O}
{\scshape J.-M. Oudom} -- {\og Th\'eor\`eme de {L}eray dans la cat\'egorie des
  alg\`ebres sur une op\'erade\fg}, \emph{C. R. Acad. Sci. Paris} \textbf{329}
  (1999), p.~101--106.

\bibitem[Pal94]{Pa}
{\scshape J.~M. Pallo} -- {\og On the listing and random generation of hybrid
  binary trees\fg}, \emph{Intern. J. Computer Math.} \textbf{50} (1994),
  p.~135--145.

\bibitem[Pir03]{Pi}
{\scshape T.~Pirashvili} -- {\og Sets with two associative operations\fg},
  \emph{C. E. J. M.} \textbf{2} (2003), p.~169--183.

\bibitem[Rem02]{Re}
{\scshape E.~Remm} -- {\og Op\'erades {L}ie-admissibles\fg}, \emph{C. R. Acad.
  Sci. Paris} \textbf{334} (2002), p.~1047--1050.

\bibitem[Ron00]{R1}
{\scshape M.~Ronco} -- {\og Primitive elements in a free dendriform
  algebra\fg}, in \emph{New trends in {H}opf algebra theory ({L}a {F}alda,
  1999)}, Contemp. Math., no. 267, Amer. Math. Soc., Providence, RI, 2000,
  p.~245--263.

\bibitem[Ron02]{R2}
\bysame , {\og Eulerian idempotents and {M}ilnor-{M}oore theorem for certain
  non-cocommutative hopf algebras\fg}, \emph{J. Algebra} \textbf{254} (2002),
  no.~1, p.~152--172.

\end{thebibliography}

\end{document}